\documentclass[11pt,a4paper,twoside]{article}
\usepackage[a4paper,inner=2cm,outer=1.8cm,top=2.4cm,bottom=2cm]{geometry}
\usepackage[T1]{fontenc}
\usepackage{lmodern}
\usepackage{amsmath,amssymb,amsthm}
\usepackage{microtype}
\usepackage{fancyhdr}
\usepackage{soul}   
\usepackage{color}  
\makeatletter
\AtBeginDocument{\let\runningauthors\@author}
\makeatother
\usepackage[colorlinks=true,linkcolor=blue,citecolor=blue,urlcolor=blue]{hyperref}
\hypersetup{pdftitle={Universal weak-type bounds for Poisson Littlewood-Paley functions on complete Riemannian manifolds},
pdfsubject={Detailed proof using sub-Markov semigroups, convex entropy, and local elliptic regularity}}
\numberwithin{equation}{section}
\newtheorem{theorem}{Theorem}[section]
\newtheorem{proposition}[theorem]{Proposition}
\newtheorem{lemma}[theorem]{Lemma}
\newtheorem{corollary}[theorem]{Corollary}
\theoremstyle{remark}
\newtheorem{remark}[theorem]{Remark}
\newcommand{\R}{\mathbb R}

\newcommand{\dd}{\,\mathrm d}
\newcommand{\norm}[2]{\lVert#1\rVert_{#2}}
\newcommand{\abs}[1]{\lvert#1\rvert}
\newcommand{\ind}{\mathbf1}
\newcommand{\Dom}{\operatorname{Dom}}
\newcommand{\supp}{\operatorname{supp}}
\newcommand{\sg}{\mathfrak g}
\newcommand{\Pmax}{P_*}
\newcommand{\Amax}{A^*}
\newcommand{\F}{\mathcal F}
\newcommand{\D}{D_{M_+}}
\allowdisplaybreaks[1]
\title{
 A geometry-independent endpoint estimate for the Littlewood--Paley $g$-function generated by the Poisson semigroup
}

\author{Bo Li, Ji Li, Jinxia Li}
\date{}
\begin{document}
\maketitle
\vspace{-1.6em}
\begin{abstract}
We prove a weak $(1,1)$ estimate for the full Poisson
Littlewood--Paley $g$-function on every complete Riemannian manifold
without boundary, with an explicit universal constant approximately equal to $2.06970$.
The estimate requires no curvature restriction, volume doubling,
heat kernel bound, or stochastic completeness.
The proof combines a precise Poisson maximal inequality with a
good-set energy estimate derived from a convex entropy identity.
The constant is sharp for the abstract combination of these two
estimates; no optimality is claimed for the operator norm.
Local elliptic regularity is used to justify the passage to
arbitrary real $L^1$ data.
As an application, we obtain weak $(1,1)$ estimates for certain
spectral multipliers represented by Poisson time derivatives with
coefficients in $L^2((0,\infty),dt/t)$.
\end{abstract}

\thanks{Keywords: Littlewood--Paley $g$-function, Riemannian manifold, spectral multiplier, entropy}

\thanks{MSC 2020: 42B25, 42B15, 58J35}

\section{Introduction and main result}

In 2003, Coulhon--Duong--Li~\cite[Theorem 1.2]{CDL}
proved that the spatial Littlewood--Paley functions associated with the heat
and Poisson semigroups are bounded on $L^p$, $1<p\le2$, on every complete
Riemannian manifold. Their weak-type $(1,1)$ version~\cite[Theorem 1.4]{CDL}, however, requires
additional geometric assumptions, including volume doubling and an
on-diagonal heat kernel upper bound, leaving the unrestricted endpoint
$p=1$ beyond the scope of their results.
A related quantitative problem concerns the optimal behavior of
Littlewood--Paley constants and dimension-free estimates. Xu~\cite{XuConstants}
addressed this problem by investigating the optimal dependence of these
constants on $p$ as $p\to1$ and $p\to\infty$.
As stated in his work,
\begin{quote}
\textit{``our main objective concerns the optimal orders of their magnitude
on $p$ as $p\to1$ and $p\to\infty$, though we also wish to have
dimension-free estimates.''}
\end{quote}
Motivated by these developments, we address the unrestricted endpoint
problem beyond the results of Coulhon--Duong--Li and investigate its
dimension-free aspect in the spirit of Xu's work. We establish a
geometry-independent weak $(1,1)$ estimate for the
Littlewood--Paley $g$-function generated by the Poisson semigroup on every complete Riemannian manifold without
boundary.
This establishes a dimension-free estimate at the weak $(1,1)$ endpoint, with a constant that is uniform over all complete Riemannian manifolds.
In particular, a warped product construction yields a dimension-free
weak $(1,1)$ bound for the first-order Poisson square function $g_1$
on the affine group (the critical case $\alpha=n-1$), thereby partially answering the
question raised by Li--Sj\"ogren \cite[Remark~3]{LiSjogren2017}; see Remark~\ref{rem:affine-group} below for more details.

We consider the following setting.
Let $M$ be a smooth, second-countable, complete Riemannian manifold
without boundary, of arbitrary finite dimension, and let $\mu$ be its
Riemannian volume measure. The differential $d$ maps scalar functions
to differential $1$-forms, whose pointwise norm is that of the
cotangent metric. The nonnegative self-adjoint Laplace--Beltrami
operator $\Delta=d^*d$ is associated with the closed Dirichlet form
$$
 q(v,w)=\int_M\langle dv,dw\rangle\dd\mu,
 \qquad \Dom q=W^{1,2}(M).
$$
In particular,
\begin{equation}\label{eq:form-domain}
 \Dom(\sqrt\Delta)=W^{1,2}(M),\qquad
 \norm{\sqrt\Delta v}{2}^2=q(v,v),
\end{equation}
and the form pairing is
\begin{equation}\label{eq:form-pairing}
 \langle\Delta v,w\rangle=q(v,w),
 \qquad v\in\Dom\Delta,\ w\in W^{1,2}(M).
\end{equation}
We use the heat semigroup $H_s=e^{-s\Delta}$ and its consistent
positive contraction extensions to $L^p(M)$, $1\le p\le\infty$.
Its strong continuity on $L^1(M)$ is standard; see~\cite{Ouhabaz,Strichartz}.
The Poisson semigroup in this paper is specifically
$P_t=e^{-t\sqrt\Delta}$ and its consistent $L^1$ extension obtained
by subordination. Both semigroups are sub-Markov:
$0\le H_s1\le1$ and $0\le P_t1\le1$. No conservation assumption
$H_s1=1$ or $P_t1=1$ is imposed. Thus $P_t$ is positive and
contractive on $L^1$ and $L^\infty$, self-adjoint on $L^2$, and
strongly continuous on $L^1$; these are consequences of the above
choice of Laplacian, not additional geometric hypotheses.
For real $f\in L^1(M)$, write $U(t,x)=P_tf(x)$, using the smooth
representative on $M_+=(0,\infty)\times M$ constructed in
Section~\ref{sec:L1}, and define the Littlewood--Paley $g$-function
\begin{equation}\label{eq:full-g}
 \sg(f)(x)=\left(\int_0^\infty
 (|t\partial_tU(t,x)|^2+|tdU(t,x)|^2)\frac{\dd t}{t}\right)^{1/2}.
\end{equation}
The integral is initially allowed to be infinite.

Our main result is the following theorem.
\begin{theorem}\label{thm:main}
Let $M$ be a complete Riemannian manifold as above.
For every real $f\in L^1(M)$ and every $\lambda>0$,
\begin{equation}\label{eq:main}
 \lambda\mu(\{x\in M:\sg(f)(x)>\lambda\})
 \le2\sqrt{C_P}\norm{f}{1},
\end{equation}
where
\begin{equation}\label{eq:exact-constants}
 C_P=\frac{\gamma(1/2,3/2)+\sqrt{3/2}\,e^{-3/2}}{\sqrt\pi},
 \qquad 2\sqrt{C_P}\approx2.0697012472,
\end{equation}
and $\gamma(\cdot, \cdot)$ denotes the lower incomplete gamma function.
In particular, $\sg(f)$ is finite almost everywhere.
The constant $2\sqrt{C_P}$ is  independent of the dimension and geometry of \(M\); no claim of optimality is made.\footnote{\ 
For a more detailed explanation, see Remark~\ref{rem:abstract-optimality} below.}
\end{theorem}

Theorem~\ref{thm:main} proves the Poisson part of
\cite[Theorem 1.4]{CDL} of Coulhon--Duong--Li without its volume doubling
and heat kernel assumptions.
In connection with the dimension-free estimates emphasized by Xu,
our contribution is an explicit universal weak $(1,1)$ bound for the scalar full Poisson $g$-function.
The uniformity extends beyond dimension: the same constant applies to all complete Riemannian manifolds, without curvature assumptions or heat kernel bounds.
The spatial and temporal Poisson functions of
\cite{CDL} are the two components of~\eqref{eq:full-g}, so each
inherits~\eqref{eq:main}. The value in~\eqref{eq:exact-constants} is
not asserted to be the optimal weak-type norm: it comes from the exact
positive subordination-weight mass $C_P$, the good-set energy bound,
and the final minimization of $a+C_P/a$ over $a>0$.

The early development of the theory includes Lusin's area integral, the work of Littlewood--Paley on Fourier series~\cite{LP},
and the later contributions of Marcinkiewicz and Zygmund; see Stein's historical account~\cite{SteinHistory}.
Stein established the classical
multidimensional theory~\cite{Stein1958} and its diffusion-semigroup
formulation~\cite{Stein}. Fefferman--Stein~\cite{FS} developed the
real-variable Hardy-space framework, while Meyer~\cite{MeyerI,Meyer}
introduced probabilistic methods for square-function inequalities.
On manifolds, curvature-based methods appear in Bakry~\cite{Bakry}
and Lohou\'e~\cite{Lohoue}; the homogeneous-type theory of
Han--Sawyer~\cite{HanSawyer} and the semigroup results of
Meda~\cite{Meda} provide further general settings.
Coulhon--Duong--Li~\cite{CDL} obtained the manifold estimates
described above. More recently, Xu developed vector-valued semigroup
theory~\cite{XuII} and sharp quantitative estimates~\cite{XuConstants}.
For the Euclidean full-gradient function, Lai~\cite{Lai} obtained
an $O(n^3)$ weak $(1,1)$ bound. Subsequent developments include
Xu's functional-calculus approach~\cite{XuFunctional},
the subordinated-semigroup maximal estimates of Bui--Cowling--Duong~\cite{BCD},
the optimal orders estimates of Lai--Li--Li~\cite{LLL},
and the obstacle decomposition for Riesz transforms in Euclidean space of Ouyang--Spector--Stockdale \cite{OSS2026}
(see the  manifold case of Chen--Jiang--Li--Li~\cite{CJLL}).

The preceding approaches use, in various combinations, Fourier and
functional calculus, real-variable decompositions, probabilistic
inequalities, and geometric or heat kernel estimates. Here we instead
organize the endpoint argument around a convex entropy and the closed
Dirichlet form.
The proof combines two classical mechanisms with precise control of
the endpoint constant. Subordination and the positive maximal ergodic
inequality give a Poisson maximal estimate; only the positive part of
the integration-by-parts weight is needed for nonnegative inputs.
For the smooth convex entropy $\Phi_R$ chosen in Section~\ref{sec:entropy}, the
harmonic extension equation and the closed form yield
$$
 F_R''(t)=\int_M\Phi_R''(U(t))(|\partial_tU(t)|^2+|dU(t)|^2)\dd\mu,
 \qquad F_R(t)=\int_M\Phi_R(U(t))\dd\mu.
$$
Our main step is to establish an entropy-based endpoint good-set distribution estimate for the full $g$-function.
We split its superlevel set into the bad part, where $P_*(|f|)>a\lambda$, and the good part, where $P_*(|f|)\le a\lambda$.
The Poisson maximal inequality controls the bad part.
The key to controlling the good part is to integrate the entropy identity above against $t\,\mathrm{d}t$, retaining the possibly nonzero entropy at infinite time.
The resulting weighted energy bound controls both the temporal and spatial derivatives and, after letting $R\downarrow a\lambda$,
yields the good-set distribution estimate: the good part has measure at most $a\|f\|_1/\lambda$.
Combining the two estimates and optimizing the maximal height gives the explicit universal constant in~\eqref{eq:exact-constants}.
Unlike the heat kernel and singular integral argument for the weak
endpoint in~\cite{CDL}, this estimate uses only positivity,
contractivity, and the form chain rule. It needs neither curvature
control nor global kernel bounds. Sub-Markovianity suffices in place
of stochastic completeness.
Local elliptic regularity ensures the smoothness of the Poisson extensions and, together with a Caccioppoli estimate, justifies the passage to $L^1$ data.
The local geometric constants involved in these arguments do not enter the final weak-type estimate.
The maximal and convex-function mechanisms have precedents
in~\cite[pp.~167--168]{Meyer}; the point here is their combination
for the full Poisson gradient with an explicit universal endpoint
constant.

As an application, let $a\in L^2((0,\infty),dt/t)$ be real-valued
and put
$$
 m_a(r)=\int_0^\infty a(t)\,tr e^{-tr}\frac{\dd t}{t}
 \quad (r>0),\qquad m_a(0)=0.
$$
We prove that
$$
 \lambda\mu(\{|m_a(\sqrt\Delta)f|>\lambda\})
 \le2\sqrt{C_P}\norm{a}{L^2(dt/t)}\norm{f}{1}
$$
for real $f\in L^1$ and $\lambda>0$. The proof uses the time derivative
in~\eqref{eq:full-g} and the Cauchy--Schwarz inequality.
The examples include differences of Poisson operators, functions of the
resolvent of $\sqrt\Delta$, and sums
$\sum_j \xi_j(P_{2^j}-P_{2^{j+1}})$ with $(\xi_j)\in\ell^2(\mathbb{Z})$.

The remainder of this paper is organized as follows.
Section~\ref{sec:maximal} establishes subordination and the Poisson
maximal inequality. Section~\ref{sec:entropy} develops spectral and
local regularity, the entropy identity and its time limits, and the
exact energy formula.
Section~\ref{sec:L1} proves the distribution estimate on $L^1\cap L^2$
and then completes the passage to arbitrary $L^1$ data.
Section~\ref{sec:multipliers} proves the spectral multiplier estimates
and gives examples.

Throughout, $L^p=L^p(M,\mu)$ and $\norm{\cdot}{p}$ denotes its norm;
the same notation is used for differential forms with their metric
norm. All scalar functions are real-valued. Inner products without
subscripts are those of $L^2$. We write $K\Subset M$ for a compact
spatial set and $\Omega\Subset M_+$ for a relatively compact open
set. On the product with metric $dt^2+g_M$, put
$\Delta_{M_+}=-\partial_t^2+\Delta$ and
$\D U=(\partial_tU,dU)$; product integrals use $\dd t\dd\mu$.
Arguments of upper-half-space functions are always ordered as
$(t,x)$. We abbreviate $U(t)=U(t,\cdot)$, and use the same convention
for $U_n(t)$ and for their spatial and time derivatives; a spatial
argument is displayed only when a pointwise statement requires it.
The maximal function is $\Pmax f=\sup_{t>0}|P_tf|$, with
the common representative specified below.

\smallskip {\noindent \textbf{AI Declaration.}}
The authors themselves formulated the problem, drawing on their
further reflections on the dimension-free weak $(1,1)$ bound for
the Riesz transform in \cite{CJLL}. OpenAI's GPT-6 assisted with
selected calculations, proof checks, language editing. 
The authors take full responsibility for the manuscript.

\section{Subordination and the Poisson maximal inequality}\label{sec:maximal}
This section derives the Poisson maximal estimate from subordination
and the positive maximal ergodic argument, keeping the kernel constant exact.

\subsection{Subordination on \texorpdfstring{$L^1$}{L1}}
The following lemma is Bochner's subordination theorem; see \cite[(1.1)]{CDL} for this Poisson kernel.
The theorem applies to the strongly continuous $L^1$ heat semigroup
and agrees with the spectral calculus on $L^1\cap L^2$.

\begin{lemma}\label{lem:subordination}
For $f\in L^1(M)$ and $t>0$, the formula
\begin{equation*}
 P_tf=\int_0^\infty\eta_t(s)H_sf\dd s,\qquad
 \eta_t(s)=\frac{t}{2\sqrt\pi}s^{-3/2}e^{-t^2/(4s)}
\end{equation*}
defines an absolutely convergent Bochner integral in $L^1$.
It agrees with the spectral Poisson semigroup on $L^1\cap L^2$.
The operators $P_t$ form a positive, sub-Markov, strongly
continuous contraction semigroup on $L^1$.
\end{lemma}

\begin{lemma}\label{lem:time-representative}
For $f\in L^1(M)$ one can choose representatives of $P_tf$
such that, outside one $\mu$-null set, $t\mapsto P_tf$
is locally absolutely continuous on $(0,\infty)$.
The maximal function $\Pmax f$ is then measurable and equals
the supremum over positive rational times.
\end{lemma}
\begin{proof}
Differentiating the scalar kernel gives
$$
 \partial_t\eta_t(s)=\eta_t(s)\left(\frac1t-\frac{t}{2s}\right).
$$
Using $r=t^2/(4s)$, we obtain
$$
 \int_0^\infty|\partial_t\eta_t(s)|\dd s
 =\frac1{t\sqrt\pi}\int_0^\infty|1-2r|r^{-1/2}e^{-r}\dd r
 =\frac1t\sqrt{\frac8{\pi e}}.
$$
For the last equality, split the integral at $r=1/2$ and use
$(2\sqrt r\,e^{-r})'=(1-2r)r^{-1/2}e^{-r}$.
Fix an integer $m\ge2$. Scalar integration over $[1/m,m]$
followed by Bochner--Fubini is legitimate because
$$
 \int_{1/m}^m\int_0^\infty
 \abs{\partial_t\eta_t(s)}\norm{H_sf}{1}\dd s\dd t
 \le2\sqrt{\frac8{\pi e}}\,\norm{f}{1}\log m<\infty.
$$
Thus on this interval
\begin{equation*}
 P_tf=P_{1/m}f+\int_{1/m}^t\partial_\tau P_\tau f\dd \tau,\qquad
 \partial_\tau P_\tau f=\int_0^\infty\partial_\tau\eta_\tau(s)H_sf\dd s
\end{equation*}
in $L^1$. Fubini applied to the last time integral defines
a representative whose scalar trajectories are absolutely
continuous for almost every $x$.
Use these representatives on the countable family $[1/m,m]$, $m\ge2$.
On overlaps the continuous representatives agree at every rational
time outside a common null set, hence at all times.
The resulting representative is locally absolutely continuous.
Its supremum is a countable supremum of measurable functions.
\end{proof}

\subsection{The Poisson maximal inequality}

For $h\ge0$ in $L^1(M)$ define
\begin{equation*}
 A_sh=\frac1s\int_0^sH_rh\dd r,\qquad \Amax h=\sup_{s>0}A_sh.
\end{equation*}
\begin{lemma}\label{lem:HDS}
For $h\ge0$ in $L^1(M)$ and $\beta>0$,
\begin{equation*}
 \beta\,\mu(\{\Amax h>\beta\})\le\norm{h}{1}.
\end{equation*}
\end{lemma}

\begin{proof}
Since $H_s$ is a positive sub-Markovian contraction
semigroup, the continuous-time Hopf maximal ergodic inequality applies
to the averages $A_sh$.
The result follows from the maximal ergodic inequality of
Dunford--Schwartz~\cite[Chapter VIII, \S7, Lemma 6, pp.~690--692]{DS}
by the standard approximation argument.
Hence, we obtain $\beta\mu(\{A^*h>\beta\})\le \|h\|_1$, for every $\beta>0$
Here the constant is $1$ rather than the factor $2$ appearing in the
general Hopf--Dunford--Schwartz maximal inequality. Indeed, in the
present setting $h\ge0$ and the semigroup is positive and
sub-Markovian, so the positive maximal argument applies directly.
The proof is complete.
\end{proof}

\begin{proposition}\label{prop:Poisson-max}
For every real $f\in L^1(M)$,
\begin{equation}\label{eq:pointwise-max}
 \Pmax f\le\Pmax(|f|)\le C_P\Amax(|f|)
 \quad\hbox{almost everywhere},
\end{equation}
where
\begin{equation}\label{eq:CP}
 C_P=\frac1{\sqrt\pi}\int_0^{3/2}
       \left(\frac32-r\right)r^{-1/2}e^{-r}\dd r
 =\frac{\gamma(1/2,3/2)+\sqrt{3/2}\,e^{-3/2}}{\sqrt\pi}.
\end{equation}
In particular,
\[
 \alpha\,\mu(\{\Pmax(|f|)>\alpha\})\le C_P\norm{f}{1},
 \qquad\alpha>0.
\]
\end{proposition}

\begin{proof}
Let $h\ge 0$ belong to $L^1(M)$. Following the integration by parts
in \cite[Section~4]{LiSjogren2017}, we rewrite the
subordination formula as
\[
 P_t h=\int_0^\infty [-s\eta_t'(s)]A_s h\dd s
 \quad\text{in }L^1,
 \qquad
 -s\eta_t'(s)=\left(\frac32-\frac{t^2}{4s}\right)\eta_t(s).
\]
For each fixed $t>0$, Fubini's theorem gives the integral representation
almost everywhere. Since $A_s h\ge0$ and $-s\eta_t'(s)\ge0$
exactly when $s\ge t^2/6$, we have
\[
P_t h = -\int_0^{t^2/6} s\eta_t'(s) A_s hds + \int_{t^2/6}^\infty [-s\eta_t'(s)] A_s hds \le C_P A^* h\quad\text{almost everywhere}.
\]
This is the exact positive-weight mass; the larger total variation
$2C_P-1$ is used only to justify absolute convergence.
Impose these inequalities first at all positive rational
times, using one null set. Lemma~\ref{lem:time-representative}
extends them to the time supremum.
Positivity gives $\abs{P_tf}\le P_t|f|$ with the same
representative convention, proving~\eqref{eq:pointwise-max}.
Applying Lemma~\ref{lem:HDS} at height $\alpha/C_P$  gives the asserted weak-type estimate, with no enlargement of $C_P$.
\end{proof}

\section{Regularity, convex entropy, and energy}\label{sec:entropy}

The analytic foundation of the proof consists of regularity of the
Poisson extension, differentiation of its convex entropy, and the
resulting exact energy formula.
Until Section~\ref{sec:L1}, fix a real
$f\in L^1(M)\cap L^2(M)$ and write $U(t)=P_tf$.

\subsection{Spectral and local regularity}
\begin{lemma}\label{lem:spectral}
For every $t>0$ and integer $k\ge0$,
$U(t)\in\Dom(\Delta^k)$. The map $U:(0,\infty)\to L^2$
is smooth and
\begin{equation}\label{eq:spectral-derivatives}
 \partial_tU(t)=-\sqrt\Delta\,U(t),\qquad \partial_t^2U(t)=\Delta U(t).
\end{equation}
The maps $t\mapsto dU(t)$ and $t\mapsto \partial_tU(t)$ are
continuous with values in their respective $L^2$ spaces.
On $M_+=(0,\infty)\times M$ the extension has a smooth
representative satisfying
\begin{equation*}
 (-\partial_t^2+\Delta)U=0.
\end{equation*}
\end{lemma}
\begin{proof}
For $t > 0$, both $\xi^k e^{-t\sqrt{\xi}}$ and $\xi^{k/2} e^{-t\sqrt{\xi}}$ are bounded for $\xi \ge 0$.
They are uniformly bounded when $t$ ranges in a compact
subinterval of $(0,\infty)$. Spectral calculus and dominated
convergence give the domain assertions and all time derivatives
in $L^2$. Equation~\eqref{eq:form-domain} gives
$$
 \norm{d(U(t)-U(s))}{2}
 =\norm{\sqrt\Delta(U(t)-U(s))}{2}\longrightarrow0,\qquad t\to s.
$$
The formula for $\partial_tU$ gives its continuity as well.
The product metric on $M_+$ has nonnegative Laplacian
$\Delta_{M_+}=-\partial_t^2+\Delta$.
Equation~\eqref{eq:spectral-derivatives}, tested against
compactly supported smooth functions on the product,
implies $\Delta_{M_+}U=0$ in distributions. The operator
$\Delta_{M_+}$ is elliptic with smooth coefficients. Local
elliptic regularity therefore makes this distribution
smooth; see~\cite[Theorem 5.1]{Hormander}.
This assertion is local and imposes no global geometric
estimate on $M$.
\end{proof}

\begin{remark}\label{rem:representatives}
The smooth product representative in Lemma~\ref{lem:spectral}
represents the spectral $P_tf$ at every $t>0$.
For a compactly supported spatial test function $\eta$,
both its pairing with the smooth representative and
$t\mapsto\langle P_tf,\eta\rangle$ are continuous.
Their equality at almost every time implies equality at
every time. The smooth representative also agrees, outside
a single spatial null set for all positive times, with
the representative of Lemma~\ref{lem:time-representative}:
use equality at rational times and continuity.
Thus the smooth extension can be used on the maximal
good set later in the proof.
\end{remark}

\subsection{The convex entropy and its properties}
Fix $A>0$. For $R>A$, let us choose a smooth even convex entropy function $\Phi_R$ such that
\begin{equation*}
 \Phi_R(0)=\Phi_R'(0)=0,\qquad
 \Phi_R''(s)=\frac1R\quad (|s|\le A),\qquad
 \Phi_R''(s)=0\quad (|s|\ge R).
\end{equation*}
It is obtained by choosing $\Phi_R''$ smoothly between $0$ and $1/R$
on $A<|s|<R$ and integrating twice from zero.
It is straightforward to verify that for every $R>A$ with $A$ fixed,
\begin{equation}\label{eq:Phi-properties}
\begin{aligned}
&0= \Phi_R(0)\leq \Phi_R(s) \leq \min \Bigl(|s|,\frac{s^2}{2R}\Bigr),\\
&|\Phi_R'(s)|\leq \min(1,|s|/R),\qquad 0\leq\Phi_R''(s)\leq\frac1R.
\end{aligned}
\end{equation}

The distinction between differentiability on all of $L^2$
and differentiability along the particular semigroup orbit
is useful. We only need the latter for the second derivative.

\begin{lemma}\label{lem:multiplier}
Let $a_j,a$ be measurable functions with
$|a_j|,|a|\le B<\infty$ and $a_j\to a$ in measure.
Then for every $w\in L^2(M)$, $ \norm{(a_j-a)w}{2}\longrightarrow0$ as $j\to\infty$.
\end{lemma}
\begin{proof}
From every subsequence one can extract a further subsequence
converging almost everywhere. On that subsequence,
$|(a_j-a)w|^2\le4B^2|w|^2$, so dominated convergence proves
the assertion.
If the full limit failed, there would be a subsequence with norms bounded away from zero, contradicting the preceding argument.
\end{proof}

\begin{proposition}\label{prop:entropy-derivatives}
For $R>A$, the function
\begin{equation*}
 F_R(t)=\int_M\Phi_R(U(t))\dd\mu
\end{equation*}
is $C^2$ on $(0,\infty)$. Its derivatives satisfy
$$ F_R'(t)=\langle\Phi_R'(U(t)),\partial_tU(t)\rangle,\qquad
 F_R''(t)=\int_M\Phi_R''(U(t))|D_{M_+}U(t)|^2\dd\mu.$$
All the displayed integrals and pairings are finite, and $F_R$ is convex.
\end{proposition}
\begin{proof}
For $R>A$, define a functional
\[
\mathcal F_R\colon L^2(M)\longrightarrow\mathbb R,\qquad
v\longrightarrow\int_M\Phi_R(v)\mathrm{d}\mu.
\]
By \eqref{eq:Phi-properties}, $\mathcal F_R$ is well defined on $L^2(M)$.
The scalar Taylor formula with $|\Phi_R''|\le1/R$ gives
\begin{equation*}
 |\F_R(v+w)-\F_R(v)-\langle\Phi_R'(v),w\rangle|
 \le\frac1{2R}\norm{w}{2}^2.
\end{equation*}
Moreover, we deduce from~\eqref{eq:Phi-properties} that
$$
 \norm{\Phi_R'(v)}{2}\le R^{-1}\norm{v}{2},\qquad
 \norm{\Phi_R'(v)-\Phi_R'(w)}{2}
 \le R^{-1}\norm{v-w}{2}.
$$
Thus $\F_R$ is continuously Fr\'echet differentiable,
with derivative represented by $\Phi_R'(v)\in L^2$.
The chain rule along the $L^2$ curve $U(t)$
proves $F_R'(t)=\langle\Phi_R'(U(t)),\partial_tU(t)\rangle$.

We next differentiate $t\mapsto\Phi_R'(U(t))$
as an $L^2$-valued map. Fix $t>0$ and, for small nonzero
$h$, put
$$
 m_h(x)=\int_0^1\Phi_R''(U(t,x)
              +\theta(U(t+h,x)-U(t,x)))\dd\theta.
$$
The pointwise fundamental theorem of calculus yields
\begin{equation}\label{eq:Phi-orbit-quotient}
 \frac{\Phi_R'(U(t+h))-\Phi_R'(U(t))}{h}
 =m_h\,\frac{U(t+h)-U(t)}h.
\end{equation}
We have $|m_h|\le1/R$.
The function $\Phi_R''$ is globally Lipschitz for fixed $R$,
since it is smooth and compactly supported.
Since $U(t+h)\to U(t)$ in $L^2$, it follows that
$m_h\to\Phi_R''(U(t))$ in measure. In fact its $L^2$
distance is bounded by a constant depending on $R$
times $\norm{U(t+h)-U(t)}{2}$.
Subtract $\Phi_R''(U(t))\partial_tU(t)$ from the right-hand side
of~\eqref{eq:Phi-orbit-quotient} and split the difference as
\begin{align*}
 m_h\left(\frac{U(t+h)-U(t)}h-\partial_tU(t)\right)
+(m_h-\Phi_R''(U(t)))\partial_tU(t).
\end{align*}
The first term tends to zero in $L^2$ by spectral
differentiability and $|m_h|\le1/R$.
The second tends to zero by Lemma~\ref{lem:multiplier}
with $w=\partial_tU(t)$. Hence
\begin{equation}\label{eq:Phi-orbit-derivative}
 \frac d{dt}\Phi_R'(U(t))
 =\Phi_R''(U(t))\partial_tU(t)\quad\hbox{in }L^2.
\end{equation}
The same multiplier lemma, together with continuity of
$U$ and $\partial_tU$ in $L^2$, shows that this derivative is
continuous. Differentiating the $L^2$ inner product in $F_R'(t)=\langle\Phi_R'(U(t)),\partial_tU(t)\rangle$ proves
\begin{equation}\label{eq:F-second-unreduced}
F_R''(t)=\int_M\Phi_R''(U(t))|\partial_tU(t)|^2\dd\mu
              +\langle\Phi_R'(U(t)),\partial_t^2U(t)\rangle
\end{equation}
and the $C^2$ assertion.
Finally, one has
$$
 \int_M\Phi_R''(U(t))|\partial_tU(t)|^2\dd\mu\le R^{-1}\norm{\partial_tU(t)}{2}^2,\qquad
 |\langle\Phi_R'(U(t)),\partial_t^2U(t)\rangle|
 \le R^{-1}\norm{U(t)}{2}\norm{\partial_t^2U(t)}{2},
$$
and the corresponding bound for $F_R'$ is immediate.
All quantities are finite by Lemma~\ref{lem:spectral}.

By Lemma~\ref{lem:spectral}, $U(t)\in\Dom\Delta\subset W^{1,2}(M)$.
Since $\Phi_R'(0)=0$ and $|\Phi_R''|\le 1/R$, the Sobolev chain
rule gives $\Phi_R'(U(t))\in W^{1,2}(M)$ and
$d(\Phi_R'(U(t)))=\Phi_R''(U(t))dU(t)$; indeed both its function and
differential have $L^2$ norms bounded by $1/R$ times the
corresponding norms of $U$.
Therefore~\eqref{eq:form-pairing} and \eqref{eq:F-second-unreduced} yield
$$
F_R''(t)
=\int_M\Phi_R''(U(t))|\partial_tU(t)|^2\dd\mu +q(U(t),\Phi_R'(U(t))) =\int_M\Phi_R''(U(t))|D_{M_+}U(t)|^2\dd\mu\ge0,
$$
which proves convexity and completes the proof.
\end{proof}

\begin{remark}
The proof does not require the integral functional $\F_R$
to be twice Fr\'echet differentiable on all of $L^2$.
Pointwise smoothness of $\Phi_R$ by itself would not
justify that assertion. Formula~\eqref{eq:Phi-orbit-derivative}
supplies exactly the second differentiation needed here,
without a global $L^\infty$ estimate for $P_tf$.
\end{remark}

\begin{lemma}\label{lem:F-endpoints}
The function $F_R$ is nonnegative, nonincreasing, and convex on
$(0,\infty)$, with $0\le F_R(t)\le\norm{f}{1}$.
Set $F_R(0)=\|\Phi_R(f)\|_1$.
Then $F_R(t)\to F_R(0)$ as $t\downarrow0$, the finite
limit $F_R(\infty)=\lim_{t\to\infty}F_R(t)$ exists, and
\begin{equation}\label{eq:F-boundaries}
 \lim_{t\downarrow0}tF_R'(t)=0,\qquad
 \lim_{t\to\infty}tF_R'(t)=0.
\end{equation}
\end{lemma}

\begin{proof}
Nonnegativity and the upper bound follow from
$0\le\Phi_R(s)\le|s|$ and $L^1$ contraction.
Convexity follows from
Proposition~\ref{prop:entropy-derivatives}.
Hence one has
\[ F_R(s)\ge F_R(t)+F_R'(t)(s-t).\]
Since $F_R$ is bounded above on $(0,\infty)$,
letting $s\to\infty$ shows that $F_R'(t)\le0$.
Thus $F_R$ is nonincreasing.

Since $\Phi_R$ is $1$-Lipschitz,
Lemma~\ref{lem:subordination} gives
\[
 |F_R(t)-F_R(0)|
 \le\norm{P_tf-f}{1}
 \longrightarrow0,\qquad t\downarrow0.
\]
Nonnegativity and monotonicity imply the existence of
$F_R(\infty)\in[0,F_R(0)]$.

Letting $s\downarrow0$ in the tangent inequality yields
\[
 0\le-tF_R'(t)\le F_R(0)-F_R(t).
\]
Taking $s=t/2$ instead gives
\[
 0\le-tF_R'(t)
 \le2(F_R(t/2)-F_R(t)).
\]
The first right-hand side tends to zero as $t\downarrow0$,
and the second tends to zero as $t\to\infty$.
This proves~\eqref{eq:F-boundaries}.
\end{proof}

\subsection{The energy formula}
\begin{proposition}\label{prop:energy}
Define the nonnegative measurable function
\begin{equation}\label{eq:J}
 J_R(x)=\int_0^\infty \Phi_R''(U(t,x))
                |tD_{M_+}U(t,x)|^2\frac{\dd t}{t}.
\end{equation}
Then
\begin{equation*}
 \int_MJ_R\dd\mu
 =F_R(0)-F_R(\infty)\le F_R(0)\le\norm{f}{1}.
\end{equation*}
\end{proposition}

\begin{proof}
The integrand in~\eqref{eq:J} is jointly measurable
and nonnegative, by the smooth product representative
of $U$.
For each integer $m\ge2$, Tonelli's theorem,
Proposition~\ref{prop:entropy-derivatives}, and
integration by parts give
\[
\int_M\int_{1/m}^{m}
   \Phi_R''(U(t,x))
   |tD_{M_+}U(t,x)|^2\frac{\dd t}{t}\dd\mu
=mF_R'(m)-m^{-1}F_R'(1/m)
          -F_R(m)+F_R(1/m).
\]
As $m\to\infty$, monotone convergence applied to
the nonnegative space-time integrand identifies
the limit of the left-hand side with
$\norm{J_R}{1}$.
Lemma~\ref{lem:F-endpoints} shows that the right-hand
side converges to $F_R(0)-F_R(\infty)$.
Consequently, one has
\[
 \int_MJ_R\dd\mu
 =F_R(0)-F_R(\infty)
 \le F_R(0)
 \le\norm{f}{1}.
\]
This also proves the finiteness of $\norm{J_R}{1}$,
without assuming it in advance.
\end{proof}

\begin{remark}
The possible nonzero value of $F_R(\infty)$ has a precise
spectral interpretation. Let $\Pi_0$ be the orthogonal
projection of $L^2$ onto $\ker\Delta$.
Spectral dominated convergence gives $P_tf\to\Pi_0f$
in $L^2$ as $t\to\infty$.
For $v,w\in L^2$, the bound $|\Phi_R'(s)|\le|s|/R$
and the scalar fundamental theorem of calculus yield
\begin{equation*}
 |\F_R(v)-\F_R(w)|
 \le\frac{\norm{v}{2}+\norm{w}{2}}{2R}\norm{v-w}{2}.
\end{equation*}
Therefore
$$
 F_R(\infty)=\int_M\Phi_R(\Pi_0f)\dd\mu.
$$
The right-hand side is finite even by the quadratic
bound~\eqref{eq:Phi-properties}.
For instance, on a compact manifold with $f\equiv1$,
$P_tf=1$ and $F_R(\infty)=\mu(M)\Phi_R(1)>0$,
whereas $J_R=0$. The infinite-time term is indispensable
for this elementary consistency check.
\end{remark}

\section{Proof of the main theorem}\label{sec:L1}
In this section, we first obtain the distribution estimate for $L^1\cap L^2$ data.
We then identify the smooth extension for arbitrary $L^1$ data,
prove local convergence of both derivatives, and pass to the limit
on the cylinders $[1/m,m]\times K$.

\subsection{The distribution estimate on \texorpdfstring{$L^1\cap L^2$}{L1 intersection L2}}
\begin{proposition}\label{prop:weak-L12}
The estimate~\eqref{eq:main} holds for real
$f\in L^1(M)\cap L^2(M)$.
\end{proposition}
\begin{proof}
Fix $\lambda,a>0$ and decompose the superlevel set into the
bad part and the good part:
$$\{\sg(f)>\lambda\}
=
\{\sg(f)>\lambda,\Pmax(|f|)>a\lambda\}
\cup \{\sg(f)>\lambda, \Pmax(|f|)\le a\lambda\}
=E_{\mathrm{bad}}\cup E_{\mathrm{good}}.
$$

On the one hand, by Proposition~\ref{prop:Poisson-max}, the bad part is controlled by
$$
\lambda\mu(E_{\mathrm{bad}})\le \lambda\mu(\{\Pmax(|f|)>a\lambda\})
\le \frac{C_P}{a}\norm{f}{1}.
$$

On the other hand, we estimate the good part by the entropy-energy argument.
Let
$$
A=a\lambda,
\qquad
\mathcal G_A=\{\Pmax(|f|)\le A\}.
$$
For $x\in\mathcal G_A$ and $t>0$, positivity and the continuous representatives
in Remark~\ref{rem:representatives} imply $|U(t,x)|\le A$, and hence $\Phi_R''(U(t,x))=1/R$ for $R>A$.
This  yields
\begin{equation*}
J_R(x)\ge\int_0^\infty \frac1R
                |tD_{M_+}U(t,x)|^2\frac{\dd t}{t}
=
\frac1R
\sg(f)(x)^2,
\qquad x\in\mathcal G_A .
\end{equation*}
Using Chebyshev's inequality and Proposition~\ref{prop:energy}, we deduce
\begin{equation*}
\lambda\mu(E_{\mathrm{good}})
\le \frac{1}{\lambda}\int_{\{\sg(f)>\lambda\}\cap\mathcal G_A}\sg(f)^2\dd\mu
\le \frac{R}{\lambda}\int_MJ_R\dd\mu
\le \frac{R}{\lambda}\norm{f}{1}.
\end{equation*}
Letting $R\downarrow A$ leads to $\lambda\mu(E_{\mathrm{good}})\le a\norm{f}{1}$.

Finally, from the estimates for the bad part and the good part, we arrive at
\begin{equation*}
\begin{aligned}
\lambda\mu(\{\sg(f)>\lambda\})\le
\lambda\mu(E_{\mathrm{bad}})+\lambda\mu(E_{\mathrm{good}})
 \le\left(\frac{C_P}{a}+a\right)\norm{f}{1}.
\end{aligned}
\end{equation*}
Optimizing in $a$ gives $a=\sqrt{C_P}$, and the corresponding coefficient is
$2\sqrt{C_P}$, as in \eqref{eq:exact-constants}.
The proof is complete.
\end{proof}

\begin{remark}
The entropy scale has been taken down to the maximal height,
which has then been optimized. The positive kernel mass
$C_P$ has been retained exactly. This does not assert that
$2\sqrt{C_P}$ is the optimal weak $(1,1)$ norm.
\end{remark}

\begin{remark}\label{rem:abstract-optimality}
The coefficient $2\sqrt{C_P}$ cannot be improved using only
the maximal weak-type bound and the good-set energy bound,
viewed as abstract inequalities.
Indeed, the two estimates used above are
\[
 A\mu(\{\Pmax(|f|)>A\})\le C_P\norm{f}{1},
 \qquad
 \int_{\mathcal G_A}\sg(f)^2\dd\mu\le A\norm{f}{1},
 \qquad A>0,
\]
where $\mathcal G_A=\{\Pmax(|f|)\le A\}$.
After normalizing $\norm{f}{1}=1$, consider the functions
$x\mapsto x$ and $x\mapsto1$ on $(0,2\sqrt{C_P})$,
with Lebesgue measure $\mu$, as abstract substitutes for
$\Pmax(|f|)$ and $\sg(f)$, respectively. For every $A>0$, we have
\[
 A\mu(\{x>A\}) =A\max\{2\sqrt{C_P}-A,0\}\le C_P,\qquad
 \int_{\{x\le A\}}1\dd\mu=\min\{A,2\sqrt{C_P}\}\le A.
\]
Thus both estimates hold, while the weak $L^1$ norm of the second function is $2\sqrt{C_P}$.
Consequently, these two estimates alone cannot yield a smaller coefficient.
No realization of this pair by a common Poisson
extension is asserted, so this example does not establish the
optimal weak $(1,1)$ norm of the Littlewood--Paley $g$-function.
\end{remark}

\subsection{The smooth representative for \texorpdfstring{$L^1$}{L1} data}
\begin{lemma}\label{lem:L1-smooth}
If $f\in L^1(M)$ is real-valued, the $L^1$ semigroup
extension $U(t)=P_tf$ admits a smooth representative on
$M_+=(0,\infty)\times M$ satisfying $\Delta_{M_+}U=0$.
For each fixed $t>0$, this smooth representative agrees
almost everywhere on $M$ with the $L^1$ element $P_tf$.
\end{lemma}
\begin{proof}
Define the amplitude truncations
\begin{equation}\label{eq:amplitude-truncation}
 f_n=\max(-n,\min(f,n)).
\end{equation}
For each $n$, we have $|f_n|\le|f|$ and $|f_n|^2\le n|f|$, so $f_n\in L^1\cap L^2$.
Moreover, $f_n\to f$ almost everywhere and $|f_n-f|\le |f|$. Hence the Lebesgue dominated convergence theorem gives $\norm{f_n-f}{1}\to0$ as $n\to \infty$.
Let $U_n(t)=P_tf_n$. Each $U_n$ is smooth and harmonic
on $M_+$ by Lemma~\ref{lem:spectral}.
For every integer $m\ge2$, contraction yields
\begin{equation}\label{eq:product-L1}
 \int_{1/m}^m\norm{U_n(t)-U(t)}{1}\dd t
 \le (m-m^{-1})\norm{f_n-f}{1}\longrightarrow0,\qquad n\to\infty.
\end{equation}
Here the strongly continuous $L^1$ curve $U(t)$ has a
jointly measurable representative, for example that of
Lemma~\ref{lem:time-representative}. It is locally integrable on $M_+$.

For every $\psi\in C_c^\infty(M_+)$, harmonicity gives
$$
 \int_{M_+} U_n\,\Delta_{M_+}\psi\dd t\dd\mu=0.
$$
Since $\Delta_{M_+}\psi$ is bounded and compactly supported,
\eqref{eq:product-L1} permits passage to the limit.
Thus $\Delta_{M_+}U=0$ in the sense of distributions.
Local elliptic regularity for the smooth product metric
makes $U$ smooth on $M_+$. This applies to distributional
solutions and does not require an a priori local $W^{1,2}$ bound.

The identification of each time slice requires one more step.
Fix $\eta\in C_c^\infty(M)$.
The pairing of $\eta$ with the smooth representative is
continuous in $t$ by smoothness on compact subsets of $M_+$.
The function $t\mapsto\langle P_tf,\eta\rangle$ is continuous
by $L^1$ strong continuity.
They agree at almost every time by Fubini and the
product-space identification, hence at every time.
At any fixed $t$ this holds for all such $\eta$, so
the two locally integrable spatial functions agree as
distributions and hence almost everywhere.
This proves the assertion for every time slice.
\end{proof}

\subsection{Local convergence of both derivatives}
\begin{lemma}\label{lem:L1-gradient-convergence}
For the truncations~\eqref{eq:amplitude-truncation}, let
$U_n(t)=P_tf_n$ and $U(t)=P_tf$ have the smooth representatives
above. If $K\Subset M$ is compact and $m\ge2$ is an integer,
then
\begin{equation}\label{eq:local-gradient-convergence}
 \int_K\int_{1/m}^m |D_{M_+}(U_n-U)|^2
 \dd t\dd\mu\longrightarrow0,\qquad n\to\infty.
\end{equation}
\end{lemma}
\begin{proof}
Apply Lemma~\ref{lem:L1-smooth} to $|f|$ and write
$V(t,x)=P_t|f|(x)$ for its nonnegative smooth representative.
Since $|f_n|\le|f|$, positivity gives $ |U_n|\le V$ and $|U|\le V$.
These inequalities first hold almost everywhere on the
product, by fixed-time positivity and Fubini.
Continuity makes them hold everywhere on $M_+$.
Put $w_n=U_n-U$.
If $\Omega\Subset M_+$ is a relatively compact open set,
then $V$ is bounded on $\overline\Omega$ and
\begin{equation}\label{eq:local-L2}
 \norm{w_n}{L^2(\Omega)}^2
 \le2\norm{V}{L^\infty(\Omega)}
           \norm{w_n}{L^1(\Omega)}\longrightarrow0,\qquad n\to\infty,
\end{equation}
by~\eqref{eq:product-L1}. This provides local $L^2$
convergence without a global $L^1\to L^2$ bound for $P_t$.

Choose $\chi\in C_c^\infty(M_+)$, equal to one on a
neighborhood of $[1/m,m]\times K$.
Each $w_n$ is smooth and harmonic.
Test its equation against the compactly supported function
$\chi^2w_n$.
With $\D w=(\partial_tw,dw)$, the product Dirichlet form and the  Caccioppoli  argument imply
\begin{align*}
 \int_{M_+}\chi^2|\D w_n|^2\dd t\dd\mu
 &=-2\int_{M_+}\chi w_n\langle\D\chi,\D w_n\rangle\dd t\dd\mu\\
 &\le\frac12\int_{M_+}\chi^2|\D w_n|^2\dd t\dd\mu
       +2\int_{M_+}|w_n|^2|\D\chi|^2\dd t\dd\mu,
\end{align*}
which implies
\begin{equation*}
 \int_{M_+}\chi^2|\D w_n|^2\dd t\dd\mu
 \le4\int_{M_+}|w_n|^2|\D\chi|^2\dd t\dd\mu\longrightarrow0,\qquad n\to\infty.
\end{equation*}
The convergence follows from~\eqref{eq:local-L2} on a
relatively compact neighborhood of $\supp\chi$.
Since \(\chi=1\) on the required cylinder, we obtain the convergence in~\eqref{eq:local-gradient-convergence}.
The proof is complete.
\end{proof}

\begin{remark}
The local bound for $V$ and the choice of $\chi$ may
depend on $M$, $f$, $K$, and $m$.
They are used only to prove convergence for a fixed
manifold and fixed cylinder. They never multiply the
weak-type bound. The limiting argument therefore does
not introduce a dimension or geometry factor
into~\eqref{eq:main}.
\end{remark}

\subsection{Truncated square functions and the final limit}
\begin{proof}[Proof of Theorem~\ref{thm:main}]
For an integer $m\ge2$, define
\begin{equation*}
 \sg_m(f)(x)=
 \left(\int_{1/m}^m
          |tD_{M_+}P_tf(x)|^2\frac{\dd t}{t}
 \right)^{1/2}.
\end{equation*}
This is finite on every compact spatial set, since the
extension is smooth on the corresponding closed cylinder.
For each fixed $x$, apply the reverse triangle inequality
in the Hilbert space
$L^2((1/m,m),t\dd t;\R\oplus T_x^*M)$.
For every compact $K\Subset M$, Lemma~\ref{lem:L1-gradient-convergence} yields
\begin{equation*}
 \norm{\sg_m(f_n)-\sg_m(f)}{L^2(K)}^2
 \le\int_K\int_{1/m}^m |t\D(U_n-U)|^2\frac{\dd t}{t}\dd\mu\longrightarrow0,\qquad n\to\infty.
\end{equation*}
Passing to a subsequence, which we do not relabel, we may assume convergence almost everywhere on $K$.
At a point where $\sg_m(f)>\lambda$,
the same strict inequality holds for all sufficiently
large indices on that subsequence. Thus
$$
 \ind_{\{\sg_m(f)>\lambda\}}
 \le\liminf_n
      \ind_{\{\sg_m(f_n)>\lambda\}}
 \quad\hbox{almost everywhere on }K.
$$
Fatou's lemma, $\sg_m(f_n)\le\sg(f_n)$,
and Proposition~\ref{prop:weak-L12} imply
$$
 \lambda\mu(K\cap\{\sg_m(f)>\lambda\})
 \le\liminf_n\lambda\mu(\{\sg_m(f_n)>\lambda\})\le2\sqrt{C_P}\lim_{n\to\infty}\norm{f_n}{1}=2\sqrt{C_P}\norm{f}{1}.
$$
The subsequence may depend on $K,m$; this
does not affect the displayed bound for the fixed limit
function. No condition on the measure of its equality
level set is necessary.

Second countability gives a countable compact exhaustion
$K_j\uparrow M$. Continuity from below of measure yields
\begin{equation*}
 \lambda\mu(\{\sg_m(f)>\lambda\})
 \le2\sqrt{C_P}\norm{f}{1}.
\end{equation*}
By nonnegativity of the integrand,
$\sg_m(f)\uparrow\sg(f)$, and the corresponding
strict superlevel sets increase to $\{\sg(f)>\lambda\}$.
Another application of continuity from below proves
\eqref{eq:main} for $f$.
For every integer $j\ge1$,
$$
 \mu(\{\sg(f)=\infty\})\le\mu(\{\sg(f)>j\})
 \le\frac{2\sqrt{C_P}}j\norm{f}{1}.
$$
Letting $j\to\infty$ proves almost everywhere finiteness.
The proof is complete.
\end{proof}

\begin{remark}\label{rem:affine-group}
Theorem~\ref{thm:main} gives a dimension-free weak $(1,1)$ bound for the
first-order Poisson square function $g_1$ on the affine group in
\cite{LiSjogren2017}.
For the critical case $\alpha=n-1$, equip
$N=\mathbb H^n\times(\R/\mathbb Z)$ with the warped metric
\[
 g_N=g_{\mathbb H^n}+y^{2(n-1)}\dd\theta^2.
\]
The projection onto $\mathbb H^n$ is $1$-Lipschitz and has compact
preimages of compact sets, so $N$ is complete. Its volume measure is
$\dd\mu_{n-1}\dd\theta$, where $\dd\mu_{n-1}=y^{-1}\dd y\dd x$ is
the right Haar measure of the affine group and $\theta$ has period one.

Circle rotations are isometries, and the Dirichlet form of $N$ on
$\theta$-independent functions is that of
$-\Delta_{\mathbb H^{(n,n-1)}}$ in \cite{LiSjogren2017}.
Thus the Poisson semigroups agree under the lift
$F(y,x,\theta)=f(y,x)$, with the identification extending to $L^1$
by contractivity. Since the lift preserves the $L^1$ norm and
$g_1(f)$ is the time component of $\sg_N(F)$,
Theorem~\ref{thm:main} gives
\[
 \lambda\,\mu_{n-1}(\{g_1(f)>\lambda\})
 \le 2\sqrt{C_P}\,\norm{f}{L^1(\mu_{n-1})},
 \qquad \lambda>0,
\]
for every real $f\in L^1(N,\mu_{n-1})$. This partially answers the
question in \cite[Remark~3]{LiSjogren2017} by settling the $g_1$
case with a constant independent of $n$.
\end{remark}

\begin{remark}
Let $M=\mathbb S^1=\mathbb R/(2\pi\mathbb Z)$ carry its standard metric and arclength measure $\dd\mu=\dd\theta$. Set
$\Delta_{\mathbb S^1}=-\partial_\theta^2$, $P_t=e^{-t\sqrt{\Delta_{\mathbb S^1}}}$, and $f(\theta)=\cos\theta$.
Since $P_tf=e^{-t}\cos\theta$, the full
Littlewood--Paley $g$-function satisfies
\[
\sg(f)(\theta)^2
=\int_0^\infty (|t\partial_tP_tf(\theta)|^2+|t\partial_\theta P_tf(\theta)|^2)\frac{\dd t}{t}
=\int_0^\infty t^2e^{-2t}\frac{\dd t}{t}
=\frac14.
\]
Thus $\|f\|_1=4$, while $\mu(\{\sg(f)>\lambda\})=2\pi$ for
$0<\lambda<1/2$ and equals zero for $\lambda\ge 1/2$. Consequently,
\[
 \sup_{\lambda>0}\lambda\mu(\{\sg(f)>\lambda\})=
\sup_{0<\lambda<1/2}2\pi\lambda
=\pi =\frac{\pi}{4}\|f\|_1.
\]
In particular, any constant $C$ for which
$\lambda\,\mu(\{\sg(f)>\lambda\})\le C\|f\|_1$
holds uniformly over the class of complete Riemannian manifolds
considered in this paper must be at least \(\pi/4\).
\end{remark}

\section{Spectral multipliers}\label{sec:multipliers}
The time derivative in~\eqref{eq:full-g} gives the following consequence of Theorem~\ref{thm:main}.
For a real function $a\in L^2((0,\infty),dt/t)$, define
$$
 m_a(r)=\int_0^\infty a(t)\,tr e^{-tr}\frac{\dd t}{t}
 \quad (r>0),\qquad m_a(0)=0.
$$
Since $\int_0^\infty(tr e^{-tr})^2\dd t/t=1/4$ for $r>0$,
the Cauchy--Schwarz inequality gives $\norm{m_a}{\infty}\le\norm{a}{L^2(dt/t)}/2$. Thus $m_a(\sqrt\Delta)$ is a bounded operator on $L^2$.

\begin{theorem}\label{thm:multiplier}
For every real $a\in L^2((0,\infty),dt/t)$, the spectral operator
$m_a(\sqrt\Delta)$ extends to real $L^1(M)$ by the formula
\begin{equation}\label{eq:multiplier-integral}
 m_a(\sqrt\Delta)f(x)=-\int_0^\infty a(t)\,t\partial_tP_tf(x)\frac{\dd t}{t}.
\end{equation}
The integral converges absolutely for almost every $x$, and
$$
 |m_a(\sqrt\Delta)f(x)|\le\norm{a}{L^2(dt/t)}\,\sg(f)(x).
$$
Consequently, for every $\lambda>0$,
\begin{equation}\label{eq:multiplier-weak}
 \lambda\mu(\{|m_a(\sqrt\Delta)f|>\lambda\})
 \le2\sqrt{C_P}\norm{a}{L^2(dt/t)}\norm{f}{1}.
\end{equation}
The constant $2\sqrt{C_P}$ is the same as in Theorem~\ref{thm:main}.
\end{theorem}
\begin{proof}
For real $f\in L^1$, the Cauchy--Schwarz inequality gives
$$
 \int_0^\infty |a(t)t\partial_tP_tf(x)|\frac{\dd t}{t}
 \le\norm{a}{L^2(dt/t)}
 \left(\int_0^\infty |t\partial_tP_tf(x)|^2
                   \frac{\dd t}{t}\right)^{1/2}
 \le\norm{a}{L^2(dt/t)}\,\sg(f)(x).
$$
By Theorem~\ref{thm:main}, the right-hand side is finite almost
everywhere. This proves absolute convergence and the pointwise bound.
The weak estimate~\eqref{eq:multiplier-weak} follows from~\eqref{eq:main}.

It remains to identify~\eqref{eq:multiplier-integral} with the spectral
operator when $f\in L^1\cap L^2$. For integers $n\ge2$, put
$a_n=a\ind_{[1/n,n]}$. The spectral theorem and
$-\partial_tP_t=\sqrt{\Delta}P_t$ give
$$
 m_{a_n}(\sqrt\Delta)f=-\int_{1/n}^n a(t)\,t\partial_tP_tf\frac{\dd t}{t}
 \quad\hbox{in }L^2.
$$
Moreover,
$$
 \norm{m_{a_n}(\sqrt\Delta)f-m_a(\sqrt\Delta)f}{2}
 \le\frac12\norm{a_n-a}{L^2(dt/t)}\norm{f}{2}
 \longrightarrow0.
$$
The pointwise integrals converge almost everywhere by absolute
convergence. Passing to an almost everywhere convergent subsequence of
the $L^2$ limit identifies the two operators.
Finally, applying~\eqref{eq:multiplier-weak} to a difference of two
$L^1$ functions proves continuity from $L^1$ to $L^{1,\infty}$.
This proves the theorem.
\end{proof}

We give several choices of $a$ and the corresponding multipliers.
First, if $0<u<v$ and $a=\ind_{[u,v]}$, then
$$
 m_a(r)=e^{-ur}-e^{-vr},\qquad
 \norm{a}{L^2(dt/t)}^2=\log(v/u).
$$
Thus $m_a(\sqrt\Delta)=P_u-P_v$.

Next, let $a(t)=t^\beta e^{-bt}$, where $\beta,b>0$.
The gamma integral gives
$$
 m_a(r)=\frac{\Gamma(\beta+1)r}{(b+r)^{\beta+1}},\qquad
 \norm{a}{L^2(dt/t)}^2=\frac{\Gamma(2\beta)}{(2b)^{2\beta}}.
$$
For example, $a(t)=te^{-t}$ gives
$m_a(\sqrt\Delta)=\sqrt{\Delta}(I+\sqrt{\Delta})^{-2}$ and $\norm{a}{L^2(dt/t)}=1/2$.

The following choice gives infinite sums with square-summable
coefficients.
\begin{corollary}\label{cor:Poisson-differences}
Let $\xi=(\xi_j)_{j\in\mathbb Z}$ be a real sequence in $\ell^2(\mathbb Z)$.
For every real $f\in L^1(M)$, the series
$$
 T_\xi f(x)=\sum_{j\in\mathbb Z} \xi_j(P_{2^j}f(x)-P_{2^{j+1}}f(x))
$$
converges absolutely almost everywhere and satisfies
\begin{equation}\label{eq:Poisson-difference-weak}
 \lambda\mu(\{|T_\xi f|>\lambda\})
 \le2\sqrt{\log2}\sqrt{C_P}\,\norm{\xi}{\ell^2}\norm{f}{1},
 \qquad\lambda>0.
\end{equation}
On $L^1\cap L^2$, it is the spectral multiplier with symbol
$$
 m_\xi(r)=\sum_{j\in\mathbb Z}
 \xi_j(e^{-2^jr}-e^{-2^{j+1}r}),\qquad r\ge0.
$$
The operator series also converges in $L^2$ for $f\in L^1\cap L^2$.
\end{corollary}
\begin{proof}
Take $a(t)=\xi_j$ on $[2^j,2^{j+1})$. Then
$$
 \norm{a}{L^2(dt/t)}^2=(\log2)\sum_{j\in\mathbb Z}|\xi_j|^2.
$$
For almost every $x$, smoothness of $P_tf$ in $t$ gives
$$
 \sum_{j\in\mathbb Z}
 |\xi_j|\,|P_{2^j}f(x)-P_{2^{j+1}}f(x)|
 \le\int_0^\infty|a(t)t\partial_tP_tf(x)|\frac{\dd t}{t}
 \le\sqrt{\log2}\,\norm{\xi}{\ell^2}\,\sg(f)(x).
$$
Hence the series converges absolutely and equals~\eqref{eq:multiplier-integral}.
Theorem~\ref{thm:multiplier} proves~\eqref{eq:Poisson-difference-weak}
and identifies its spectral symbol.
Approximating $\xi$ by finitely supported sequences in $\ell^2$ and using the estimate $\norm{m_a}{\infty}\le\norm{a}{L^2(dt/t)}/2$ proves convergence in $L^2$.
The proof is complete.
\end{proof}

For a concrete sequence in $\ell^2\setminus\ell^1$, take
$\xi_j=(-1)^j(j+1)^{-\alpha}$ for $j\ge1$ and $\xi_j=0$ for $j\le0$,
where $1/2<\alpha\le1$. The corresponding coefficient is
$$
a(t)=
\begin{cases}
\dfrac{(-1)^j}{(j+1)^\alpha},\quad  &\text{if }2^j\le t<2^{j+1},\ j\ge1,\\
0,&\text{if }0<t<2.
\end{cases}
$$
It belongs to $L^2(dt/t)$, but not to $L^1(dt/t)$.
Thus~\eqref{eq:Poisson-difference-weak} applies even when
$\|\xi\|_{\ell^1}=\infty$.

We conclude this paper with the following remark.

\begin{remark}\label{rem:sigma-Poisson}
The Littlewood--Paley $g$-function estimates and the associated
spectral multiplier estimates established in this paper extend
to $\sigma$-Poisson operators. We outline the main modifications,
omitting the details.

For $0<\sigma<1$, the Caffarelli--Silvestre extension problem
\cite{CaffarelliSilvestre}, in the operator formulation of
Stinga--Torrea \cite{StingaTorrea}, gives the $\sigma$-Poisson
extension $U(t)=P_t^\sigma f$, defined by
\[
 P_t^\sigma f
 =\frac{t^{2\sigma}}{4^\sigma\Gamma(\sigma)}
   \int_0^\infty s^{-1-\sigma}
   e^{-t^2/(4s)}H_sf\dd s.
\]
This extension satisfies
\[
 \partial_t^2U(t)
 +\frac{1-2\sigma}{t}\partial_tU(t)
 -\Delta U(t)=0.
\]
Substituting this equation into~\eqref{eq:F-second-unreduced}
and applying the form identity and the Sobolev chain rule
yields the entropy identity
\[
 F_R''(t)+\frac{1-2\sigma}{t}F_R'(t)
 =\int_M\Phi_R''(U(t))
 |D_{M_+}U(t)|^2\dd\mu.
\]
The Caccioppoli argument in the proof of
Lemma~\ref{lem:L1-gradient-convergence} is carried out
with respect to the weighted measure
$t^{1-2\sigma}\dd t\dd\mu$; see also the related cutoff
argument in the proof of \cite[Lemma~3.2]{LLLMS2024}.
For the spectral multiplier estimates, replace $tr e^{-tr}$
in the definition of $m_a$ by $\psi_\sigma(tr)$, where
\[
 \psi_\sigma(r)
 =\frac{r^2}{2\Gamma(\sigma)}
   \int_0^\infty s^{\sigma-2}
   e^{-s-r^2/(4s)}\dd s
 \quad(r>0),
 \qquad
 \psi_\sigma(0)=0.
\]
With the corresponding adjustments to the remaining arguments
and formulas, one obtains estimates whose constants may depend
on $\sigma$ but remain independent of the dimension and
geometry of $M$.
\end{remark}

\smallskip

{\footnotesize{\noindent \textbf{Acknowledgments.}
We thank Prof. Hong-Quan Li for pointing us to \cite{LiSjogren2017}, which allowed us to simplify the proof of Proposition~\ref{prop:Poisson-max}.
Bo Li is supported by NNSF of China (12471094, 12571121), Zhejiang NSF of China (LMS26A010015).
Ji Li is supported by ARC (DP260100485).
Jinxia Li is supported by NNSF of China (12101199, 12571101), Henan NSF of China (232300421142, 242300421387).
}}

\begingroup
\footnotesize
\setlength{\parskip}{0pt}
\setlength{\itemsep}{0pt}
\setlength{\parsep}{0pt}
\setlength{\topsep}{0pt}
\setlength{\partopsep}{0pt}
\setlength{\labelsep}{0.3em}
\renewcommand{\baselinestretch}{0.95}\selectfont
\emergencystretch=1em

\endgroup

\footnotesize{
\noindent Bo Li (bli@zjxu.edu.cn), Institute of Mathematics, Jiaxing University, Jiaxing 314001, China.

\smallskip

\noindent Ji Li (ji.li@mq.edu.au), School of Mathematical and Physical Sciences, Macquarie University, NSW 2109, Australia.

\smallskip

\noindent Jinxia Li (jinxiali@hpu.edu.cn), School of Mathematics and Information Science, Henan Polytechnic University, Jiaozuo 454003, China.
}

\end{document}